%
\def\autori{P.\ DALL'AGLIO and G.\ DAL MASO}
\def\titolo{Some properties of the solutions of obstacle problems with
measure data} 




\font\sixrm=cmr6
\newcount\tagno \tagno=0                        
\newcount\thmno \thmno=0                        
\newcount\bibno \bibno=0                        
\newcount\chapno\chapno=0                       
\newcount\verno            
\newif\ifproofmode
\proofmodetrue
\newif\ifwanted
\wantedfalse
\newif\ifindexed
\indexedfalse

\def\ifundefined#1{\expandafter\ifx\csname+#1\endcsname\relax}

\def\Wanted#1{\ifundefined{#1} \wantedtrue 
\immediate\write0{Wanted 
#1 
\the\chapno.\the\thmno}\fi}

\def\Increase#1{{\global\advance#1 by 1}}

\def\Assign#1#2{\immediate
\write1{\noexpand\expandafter\noexpand\def
 \noexpand\csname+#1\endcsname{#2}}\relax
 \global\expandafter\edef\csname+#1\endcsname{#2}}

\def\pAssign#1#2{\write1{\noexpand\expandafter\noexpand\def
 \noexpand\csname+#1\endcsname{#2}}}

\def\lPut#1{\ifproofmode\llap{\hbox{\sixrm #1\ \ \ }}\fi}
\def\rPut#1{\ifproofmode$^{\hbox{\sixrm #1}}$\fi}



\def\chp#1{\global\tagno=0\global\thmno=0\Increase\chapno
\Assign{#1}
{\the\chapno}{\lPut{#1}\the\chapno}}


\def\thm#1{\Increase\thmno 
\Assign{#1}
{\the\chapno.\the\thmno}\the\chapno.\the\thmno\rPut{#1}}


\def\frm#1{\Increase\tagno
\Assign{#1}{\the\chapno.\the\tagno}\lPut{#1}
{\the\chapno.\the\tagno}}


\def\bib#1{\Increase\bibno
\Assign{#1}{\the\bibno}\lPut{#1}{\the\bibno}}


\def\pgp#1{\pAssign{#1/}{\the\pageno}}


\def\ix#1#2#3{\pAssign{#2}{\the\pageno}
\immediate\write#1{\noexpand\idxitem{#3}
{\noexpand\csname+#2\endcsname}}}

 
\def\rf#1{\Wanted{#1}\csname+#1\endcsname\relax\rPut {#1}}


\def\rfp#1{\Wanted{#1}\csname+#1/\endcsname\relax\rPut{#1}}

\input \jobname.auxi
\Increase\verno
\immediate\openout1=\jobname.auxi

\immediate\write1{\noexpand\verno=\the\verno}

\ifindexed
\immediate\openout2=\jobname.idx
\immediate\openout3=\jobname.sym 
\fi


\font\twelverm=cmr12
\font\twelvei=cmmi12
\font\twelvesy=cmsy10
\font\twelvebf=cmbx12
\font\twelvett=cmtt12
\font\twelveit=cmti12
\font\twelvesl=cmsl12

\font\ninerm=cmr9
\font\ninei=cmmi9
\font\ninesy=cmsy9
\font\ninebf=cmbx9
\font\ninett=cmtt9
\font\nineit=cmti9
\font\ninesl=cmsl9

\font\eightrm=cmr8
\font\eighti=cmmi8
\font\eightsy=cmsy8
\font\eightbf=cmbx8
\font\eighttt=cmtt8
\font\eightit=cmti8
\font\eightsl=cmsl8

\font\sixrm=cmr6
\font\sixi=cmmi6
\font\sixsy=cmsy6
\font\sixbf=cmbx6

\catcode`@=11 
\newskip\ttglue

\def\twelvepoint{\def\rm{\fam0\twelverm}
\textfont0=\twelverm  \scriptfont0=\ninerm  
\scriptscriptfont0=\sevenrm
\textfont1=\twelvei  \scriptfont1=\ninei  \scriptscriptfont1=\seveni
\textfont2=\twelvesy  \scriptfont2=\ninesy  
\scriptscriptfont2=\sevensy
\textfont3=\tenex  \scriptfont3=\tenex  \scriptscriptfont3=\tenex
\textfont\itfam=\twelveit  \def\it{\fam\itfam\twelveit}%
\textfont\slfam=\twelvesl  \def\sl{\fam\slfam\twelvesl}%
\textfont\ttfam=\twelvett  \def\tt{\fam\ttfam\twelvett}%
\textfont\bffam=\twelvebf  \scriptfont\bffam=\ninebf
\scriptscriptfont\bffam=\sevenbf  \def\bf{\fam\bffam\twelvebf}%
\tt  \ttglue=.5em plus.25em minus.15em
\normalbaselineskip=15pt
\setbox\strutbox=\hbox{\vrule height10pt depth5pt width0pt}%
\let\sc=\tenrm  \let\big=\twelvebig  \normalbaselines\rm}

\def\tenpoint{\def\rm{\fam0\tenrm}
\textfont0=\tenrm  \scriptfont0=\sevenrm  \scriptscriptfont0=\fiverm
\textfont1=\teni  \scriptfont1=\seveni  \scriptscriptfont1=\fivei
\textfont2=\tensy  \scriptfont2=\sevensy  \scriptscriptfont2=\fivesy
\textfont3=\tenex  \scriptfont3=\tenex  \scriptscriptfont3=\tenex
\textfont\itfam=\tenit  \def\it{\fam\itfam\tenit}%
\textfont\slfam=\tensl  \def\sl{\fam\slfam\tensl}%
\textfont\ttfam=\tentt  \def\tt{\fam\ttfam\tentt}%
\textfont\bffam=\tenbf  \scriptfont\bffam=\sevenbf
\scriptscriptfont\bffam=\fivebf  \def\bf{\fam\bffam\tenbf}%
\tt  \ttglue=.5em plus.25em minus.15em
\normalbaselineskip=12pt
\setbox\strutbox=\hbox{\vrule height8.5pt depth3.5pt width0pt}%
\let\sc=\eightrm  \let\big=\tenbig  \normalbaselines\rm}

\def\ninepoint{\def\rm{\fam0\ninerm}
\textfont0=\ninerm  \scriptfont0=\sixrm  \scriptscriptfont0=\fiverm
\textfont1=\ninei  \scriptfont1=\sixi  \scriptscriptfont1=\fivei
\textfont2=\ninesy  \scriptfont2=\sixsy  \scriptscriptfont2=\fivesy
\textfont3=\tenex  \scriptfont3=\tenex  \scriptscriptfont3=\tenex
\textfont\itfam=\nineit  \def\it{\fam\itfam\nineit}%
\textfont\slfam=\ninesl  \def\sl{\fam\slfam\ninesl}%
\textfont\ttfam=\ninett  \def\tt{\fam\ttfam\ninett}%
\textfont\bffam=\ninebf  \scriptfont\bffam=\sixbf
\scriptscriptfont\bffam=\fivebf  \def\bf{\fam\bffam\ninebf}%
\tt  \ttglue=.5em plus.25em minus.15em
\normalbaselineskip=11pt
\setbox\strutbox=\hbox{\vrule height8pt depth3pt width0pt}%
\let\sc=\sevenrm  \let\big=\ninebig  \normalbaselines\rm}

\def\eightpoint{\def\rm{\fam0\eightrm}
\textfont0=\eightrm  \scriptfont0=\sixrm  \scriptscriptfont0=\fiverm
\textfont1=\eighti  \scriptfont1=\sixi  \scriptscriptfont1=\fivei
\textfont2=\eightsy  \scriptfont2=\sixsy  \scriptscriptfont2=\fivesy
\textfont3=\tenex  \scriptfont3=\tenex  \scriptscriptfont3=\tenex
\textfont\itfam=\eightit  \def\it{\fam\itfam\eightit}%
\textfont\slfam=\eightsl  \def\sl{\fam\slfam\eightsl}%
\textfont\ttfam=\eighttt  \def\tt{\fam\ttfam\eighttt}%
\textfont\bffam=\eightbf  \scriptfont\bffam=\sixbf
\scriptscriptfont\bffam=\fivebf  \def\bf{\fam\bffam\eightbf}%
\tt  \ttglue=.5em plus.25em minus.15em
\normalbaselineskip=9pt
\setbox\strutbox=\hbox{\vrule height7pt depth2pt width0pt}%
\let\sc=\sixrm  \let\big=\eightbig  \normalbaselines\rm}

\def\twelvebig#1{{\hbox{$\textfont0=\twelverm\textfont2=\twelvesy
	\left#1\vbox to10pt{}\right.\n@space$}}}
\def\tenbig#1{{\hbox{$\left#1\vbox to8.5pt{}\right.\n@space$}}}
\def\ninebig#1{{\hbox{$\textfont0=\tenrm\textfont2=\tensy
	\left#1\vbox to7.25pt{}\right.\n@space$}}}
\def\eightbig#1{{\hbox{$\textfont0=\ninerm\textfont2=\ninesy
	\left#1\vbox to6.5pt{}\right.\n@space$}}}
 
\def\displayliness#1{\null\,\vcenter{\openup1\jot \m@th
  \ialign{\strut\hfil$\displaystyle{##}$\hfil
    \crcr#1\crcr}}\,}
	       
\def\displaylinesno#1{\displ@y \tabskip=\centering
   \halign to\displaywidth{ \hfil$\@lign \displaystyle{##}$ \hfil
	\tabskip=\centering
     &\llap{$\@lign##$}\tabskip=0pt \crcr#1\crcr}}
		     
\def\ldisplaylinesno#1{\displ@y \tabskip=\centering
   \halign to\displaywidth{ \hfil$\@lign \displaystyle{##}$\hfil
	\tabskip=\centering
     &\kern-\displaywidth
     \rlap{$\@lign##$}\tabskip=\displaywidth \crcr#1\crcr}}

\catcode`@=12 

\def\parag#1#2{\goodbreak\bigskip\bigskip\noindent
		   {\bf #1.\ \ #2}
		   \nobreak\bigskip} 
\def\intro#1{\goodbreak\bigskip\bigskip\goodbreak\noindent
		   {\bf #1}\nobreak\bigskip\nobreak}
\long\def\th#1#2{\goodbreak\bigskip\noindent
		{\bf Theorem #1.\ \ \it #2}}
\long\def\lemma#1#2{\goodbreak\bigskip\noindent
		{\bf Lemma #1.\ \ \it #2}}

\long\def\cor#1#2{\goodbreak\bigskip\noindent
		{\bf Corollary #1.\ \ \it #2}}
\long\def\defin#1#2{\goodbreak\bigskip\noindent
		  {\bf Definition #1.\ \ \rm #2}}
\long\def\rem#1#2{\goodbreak\bigskip\noindent
		 {\bf Remark #1.\ \ \rm #2}}
\long\def\ex#1#2{\goodbreak\bigskip\noindent
		 {\bf Example #1.\ \ \rm #2}}

\def\proof{\vskip.4cm\noindent{\it Proof.\ \ }}

\def\sqr#1#2{\vbox{
   \hrule height .#2pt 
   \hbox{\vrule width .#2pt height #1pt \kern #1pt 
      \vrule width .#2pt}
   \hrule height .#2pt }}
\def\square{\sqr74}

\def\endproof{{\unskip\nobreak\hfill \penalty50
\hskip2em\hbox{}\nobreak\hfill $\square$ \goodbreak
\parfillskip=0pt  \finalhyphendemerits=0}}

\mathchardef\emptyset="001F
\mathchardef\hyphen="002D


\def\rightheadline{\eightpoint\hfil\titolo
\hfil\tenrm\folio} 
\def\leftheadline{\tenrm\folio\hfil\eightpoint
\autori \hfil}
\def\zeroheadline{\hfill} 
%
\headline={\ifnum\pageno=0 \zeroheadline
\else\ifodd\pageno\rightheadline
\else\leftheadline\fi\fi}

\nopagenumbers
\magnification=1200
\baselineskip=15pt
\topskip=15pt 
\hfuzz=2pt
\parindent=2em
\mathsurround=1pt
\tolerance=1000

\pageno=0
\hsize 14truecm
\vsize 25truecm
\hoffset=0.8truecm
\voffset=-1.55truecm

\null
\vskip 2.8truecm
{\twelvepoint
\baselineskip=1.7\baselineskip
\centerline{\bf SOME PROPERTIES OF THE SOLUTIONS}
\centerline{\bf  OF OBSTACLE PROBLEMS}
\centerline{\bf WITH MEASURE DATA}
}
\vskip2truecm
\centerline{Paolo DALL'AGLIO}
\medskip
\centerline{Gianni DAL MASO}

\vfil

{\eightpoint
\baselineskip=1.2\baselineskip
\centerline{\bf Abstract}
\bigskip
\noindent
We study some properties of the obstacle reactions associated with the
solutions of unilateral obstacle problems with measure data. These
results allow us to prove that, under very weak
assumptions on the obstacles, the solutions do not depend on the
components of the negative parts of the data which are concentrated on
sets of capacity zero. The proof is based on a careful analysis of 
the behaviour of the potentials of two mutually singular measures 
near the points where both potentials tend to infinity.
\par
}
\vfil
\medskip
\centerline{SISSA, via Beirut 4, 34014 Trieste, Italy}
\medskip
\centerline{e-mail:\ \thinspace{\tt aglio@sissa.it}\thinspace,\ 
\ {\tt dalmaso@sissa.it}}
\vskip 1truecm
\centerline {Ref. S.I.S.S.A.
 126/98/M (November 98)}
\vskip 1truecm
\eject
%
%
%
\pageno=1
\topskip=25pt
\vsize 22.5truecm
\hsize 16.2truecm
\hoffset=0truecm
\voffset=0.5truecm
\baselineskip=15pt
\hfuzz=2pt
\parindent=2em
\mathsurround=1pt

\def\R{{\rm I\! R}}
\def\Rn{{\rm I\! R}^N}
\def\Luno{{\rm L}^1(\Omega)}
\def\Linf{{\rm L}^\infty(\Omega)}
\def\Hunozero{{\rm H}^1_0(\Omega)}   
\def\Huno{{\rm H}^1(\Omega)}
\def\Hduale{{\rm H}^{\hbox{\kern 1truept \rm -}\kern -1truept1}\kern
	-1truept(\Omega)}
\def\Mb{{\cal M}_b(\Omega)}
\def\Mbp{{\cal M}^+_b(\Omega)}
\def\Mbo{{\cal M}^0_b(\Omega)}
\def\Mbop{{\cal M}^{0,+}_b(\Omega)}
\def\Inters{\Mb\cap\Hduale}
\def\K{{\cal K}_\psi(\Omega)}
\def\intl{\int\limits}
\def\dualita#1#2{\langle#1,#2\rangle}
\def\meas{{\rm meas}}
\def\um{u_{\mu}}
\def\un{u_{\nu}}
\def\ul{u_{\lambda}}

\def\Gao{G^{\cal A}_\Omega}

\def\qe{\,\hbox{\ q.e.\ in }\,\Omega}
\def\e{\varepsilon}
\def\ave{\mathop{-\hskip-.38cm\intop}\limits}
\def\cc{\subset\subset}

\def\crr{\cr\noalign{\medskip}}

\proofmodefalse 

\parag{\chp{intro}}{Introduction}

Given a regular bounded open set $\Omega$ of ${\bf R}^N$, $N\ge2$, and a 
linear elliptic operator ${\cal A}$ of the form
$$
{\cal A}u =-\sum_{i,j=1}^N D_i ( a_{ij} D_j u )\,,
\eqno(\frm{ellop})
$$
with $a_{ij}\in \Linf$, we study some properties of 
the solution of the obstacle problem for the operator ${\cal A}$ in
$\Omega$ with 
homogeneous Dirichlet boundary conditions on $\partial\Omega$, when the datum 
$\mu$ is a bounded Radon measure on $\Omega$ and the obstacle 
$\psi$ is an arbitrary function on $\Omega$. According to 
[\rf{DAL-LEO}], a function $u$ is a solution of this  problem, which will 
be denoted by $OP(\mu,\psi)$, if $u$ is the smallest function with the 
following properties: $u\ge\psi$ in $\Omega$ and $u$ is a solution 
in the sense of Stampacchia [\rf{STA}] of a problem of the form
$$
\cases{{\cal A}u=\mu+\lambda&in $\Omega\,$,\crr
u=0&on $\partial\Omega\,$,
\cr}
\eqno(\frm{DM1})
$$
for some bounded Radon measure $\lambda\ge0$. 
The measure $\lambda$ which corresponds to the solution of the 
obstacle problem is called the obstacle reaction.

Existence and uniqueness of the solution of $OP(\mu,\psi)$ have 
been proved in [\rf{DAL-LEO}], provided that there exists a measure
$\lambda$ such that the solution of (\rf{DM1}) is greater than or 
equal to $\psi$.
These results have been extended to the non-linear 
case in~[\rf{LEO}], when $\mu$ vanishes on all sets with capacity zero.
For a different approach to obstacle problems for non-linear 
operators with measure data
see [\rf{BOC-GAL}], [\rf{BOC-CIR1}], 
[\rf{BOC-CIR2}], [\rf{OPP-ROS1}], and [\rf{OPP-ROS2}].

If the measure $\mu$ belongs to the dual $\Hduale$ of the 
Sobolev space $\Hunozero$, and
if there exists a function $w\in \Hunozero$ above 
the obstacle $\psi$, then  the solution of the obstacle 
problem $OP(\mu,\psi)$ according to the previous definition coincides 
with the solution $u$ of the variational inequality
$$
\cases{u\in \Hunozero\,,\ u\ge\psi\,,&\crr
\langle{\cal A}u, v-u\rangle \ge \langle\mu, v-u\rangle &\crr
\forall v\in \Hunozero\,,\ v\ge\psi\,,&\cr}
\eqno(\frm{DM2})
$$
where $\langle\cdot,\cdot\rangle$ denotes the duality pairing between
$\Hduale$ and $\Hunozero$.
In this case the obstacle reaction $\lambda$ belongs to 
$\Hduale$. It is concentrated on the contact set 
$\{u=\psi\}$ if $\psi$ is continuous, or, more in general, quasi 
upper semicontinuous.

An important role in this problem is played by the space 
${\cal M}_b^0(\Omega)$ of all bounded Radon measures on $\Omega$ 
which are absolutely continuous with respect to the harmonic capacity.
If the datum $\mu$ belongs to ${\cal M}_b^0(\Omega)$, so does the 
obstacle reaction, provided that there exists a measure 
$\lambda\in {\cal M}_b^0(\Omega)$ such that the solution of (\rf{DM1}) 
is greater than or equal to $\psi$ (see [\rf{DAL-LEO}], Theorem~7.5). 
In this case the obstacle reaction is concentrated on the contact set 
$\{u=\psi\}$, whenever the obstacle $\psi$ is quasi upper semicontinuous
(see [\rf{LEO}], Theorem~2.9). Example~\rf{delta}, which is 
a variant of an example proposed by L.~Orsina and A.~Prignet, shows that
this is not 
always true when $\mu$ is not absolutely continuous with respect to 
the harmonic capacity.

Using the linearity of the operator ${\cal A}$, it is easy to see that the 
obstacle reaction belongs to ${\cal M}_b^0(\Omega)$ and is concentrated 
on the contact set $\{u=\psi\}$, whenever $\psi$ is quasi upper 
semicontinuous and just the negative part $\mu^-$ of 
$\mu$ belongs to ${\cal M}_b^0(\Omega)$. Therefore we concentrate our 
attention on the 
case $\mu^-\notin{\cal M}_b^0(\Omega)$. Then $\mu^-$ can be 
decomposed as $\mu^-=\mu^-_a+\mu^-_s$, where 
$\mu^-_a \in{\cal M}_b^0(\Omega)$ and $\mu^-_s$ is concentrated on a 
set 
of capacity zero. We assume that the obstacle $\psi$ satisfies the 
estimates $-v-w \le \psi\le v$,
where $w\in {\rm H}^1(\Omega)$ and $v$ is the solution in the 
sense 
of Stampacchia of a problem of the form
$$
\cases{{\cal A}v=\nu &in $\Omega\,$,\crr
v=0&on $\partial\Omega\,$,
\cr}
\eqno(\frm{DM4})
$$
with $\nu\in{\cal M}_b^0(\Omega)$. We prove 
(Theorem~\rf{lostesso}) that 
the obstacle problems $OP(\mu,\psi)$ and $OP(\mu^+-\mu^-_a,\psi)$ 
have the same solution $u$, while the corresponding obstacle reactions 
$\lambda$ and $\lambda_0$ satisfy $\lambda=\lambda_0+\mu^-_s$. 
This shows that, under these assumptions, the solution $u$ of  
$OP(\mu,\psi)$ does not depend on $\mu^-_s$, while the obstacle 
reaction has the form $\lambda_0+\mu^-_s$, where $\lambda_0$ is a 
non-negative measure in ${\cal M}_b^0(\Omega)$. This measure is 
concentrated on the contact set $\{u=\psi\}$ whenever
the obstacle $\psi$ is 
quasi upper semicontinuous (Theorem~\rf{complementarieta}). 

These results will be used in a forthcoming paper [\rf{DAL}] to study the 
dependence of the solutions on the obstacles. Their proof relies on a 
variant (Lemma~\rf{unico}) of the following result, which has an 
intrinsic interest. Let 
$u_\mu$ and $u_\nu$ be the solutions of (\rf{DM4}) 
corresponding to the measures $\mu$ and $\nu$, which are not assumed to
belong to ${\cal M}_b^0(\Omega)$. Suppose that $\mu^+\!\perp\nu$ and 
$u_\mu\le u_\nu$. Then $\mu^+\!\in {\cal M}_b^0(\Omega)$. This result is
obtained by 
investigating the behaviour of the potentials of two mutually singular 
measures near their singular points (Lemmas~\rf{rapporto} 
and~\rf{rapportoA}).

\vfill
\eject

\parag{\chp{fattinoti}}{Notation and preliminary results}

Let us fix a bounded open set $\Omega$ 
in $\Rn$, $N\geq 2$. We assume that $\Omega$ satisfies the
following regularity condition, considered by Stampacchia in
[\rf{STA}]: there exists a constant $\alpha>0$ such that
$$
\meas(B_{r}(x)\setminus \Omega) \ge \alpha \,\meas(B_{r}(x))\,,
$$
for every $x\in\partial\Omega$ and for every $r>0$, where $B_{r}(x)$ 
denotes the open ball with centre $x$ and radius~$r$.

Let ${\cal A}$ be the linear elliptic 
operator introduced in (\rf{ellop}), 
where $(a_{ij})$ is an $N{\times}N$ matrix
of functions in $L^\infty(\Omega)$, and, for a suitable
constant $\beta>0$,
$$
\sum_{i,j=1}^{N} a_{ij}(x)\xi_i\xi_j\geq\beta\, |\xi|^2,
$$
for a.e.\ $x\in\Omega$ and for every $\xi\in\R^N$.

In order to include in our analysis also the case of thin
obstacles, it is convenient to introduce the notions of capacity and 
of quasi continuous representative of a Sobolev function.
Given a set $E\subseteq\Omega$, its
{\it capacity\/} with respect to $\Omega$ is defined by
$$
{\rm cap}(E)=\inf\intl_{\Omega} |\nabla v|^2 dx \,,
$$
where $v$ runs over all functions $v\in\Hunozero$ such that $v\geq 1$ 
a.e.\ in a neighbourhood of $E$.
We say that a property holds {\it quasi everywhere\/} 
(abbreviated as {\it q.e.\/}) when it holds everywhere except
on a set of capacity zero.
A function $v\colon \Omega\to{\overline{\R}}$ is {\it quasi 
continuous\/} 
(resp.\ {\it quasi upper semicontinuous\/}) if, 
for every $\e>0$, there exists a set
$E\subseteq\Omega$, with ${\rm cap}(E)<\e$, such that 
$v|_{\Omega\setminus E}$ is continuous
(resp.\ upper semicontinuous) in
$\Omega\setminus E$.
We recall also that, if $u$ and $v$ are quasi continuous functions and
$u\leq v$ a.e.\ in $\Omega$, then $u\leq v$ q.e.\ in $\Omega$.

Every function $u\in\Hunozero$ has a {\it quasi continuous 
representative\/},
i.e., a quasi continuous function $\tilde u$ which is equal to $u$ a.e.\ 
in $\Omega$. We
shall always identify $u$ with its quasi continuous representative $\tilde u$,
which is uniquely defined quasi everywhere in~$\Omega$.
A self-contained presentation of all these notions can be found, 
for instance, in Chapters~4 of [\rf{EVA-GAR}] and [\rf{HEI-KIL}].

Let us fix a function $\psi\colon\Omega\to{\overline{\R}}$, and the 
corresponding convex set
$$
{\cal K}_{\psi}(\Omega):=\{z\ \hbox{quasi continuous in } \Omega : 
z\geq \psi {\hbox{ q.e.\ in }}\Omega\}\,. 
$$

In their natural setting, obstacle problems are part of
the theory of Variational
Inequalities (for which we refer to the books [\rf{BAI-CAP}],
[\rf{KIN-STA}], and [\rf{TRO}]).
For any $\mu\in\Hduale$ the
variational inequality with obstacle $\psi$
$$
\cases{u\in\K\cap\Hunozero\,,\crr
                      \dualita{{\cal A}u}{v-u}
                      \geq\dualita{\mu}{v-u} \crr
		\forall v\in\K\cap\Hunozero\,,
	\cr}
\eqno(\frm{vi})
$$
which will be indicated by $VI(\mu,\psi)$,
has a unique solution $u$, whenever the set $\K\cap \Hunozero$ is nonempty,
i.e.,
$$
\hbox{there exists }\, w\in\Hunozero\, \hbox{ such that } \,
w\geq\psi\qe\,.\eqno(\frm{zeta})
$$
In this case we say that the obstacle is {\it $VI$-admissible\/}.

Among all classical results, we recall
that the solution of $VI(\mu,\psi)$ is also characterized as
the smallest function $u\in\Hunozero$ such that
$$
\cases{
{\cal A}u-\mu\geq 0\,
	\hbox{ in }\,{\cal D}'(\Omega)\,,\crr
u\geq\psi\qe\,.\cr}\eqno(\frm{caratt})
$$
Then $\lambda:={\cal A}u-\mu$ is a non-negative measure, that is called
the {\it obstacle reaction\/} associated with~$u$.
\bigskip

Let $\Mb$ be the space of all bounded Radon measures on $\Omega$, 
and let $\Mbo$ be
the subspace of all measures of $\Mb$ which vanish on all sets 
of capacity zero. The corresponding cones of non-negative measures
will be denoted by $\Mbp$ and $\Mbop$, respectively. 
Recall that $\Hduale\not\subseteq\Mb$, but 
$\Hduale\cap\Mb\subseteq\Mbo$.
Any measure $\mu\in\Mb$ can be decomposed as $\mu=\mu_a+\mu_s$, 
where $\mu_a\in\Mbo$ and $\mu_s$ is concentrated on a set of capacity
zero (see [\rf{FUK-SAT-TAN}]).

\bigskip
When the datum is a measure, equations and inequalities can not be studied
in the variational framework, and the usual notion of solution in the 
sense of distributions does not guarantee uniqueness when the 
coefficients are discontinuous, as shown by a 
celebrated counterexample due to J.~Serrin [\rf{SER}].
To overcome these difficulties, Stampacchia introduced
in~[\rf{STA}] the following notion of solution, obtained by duality.

\defin{\thm{stamp}}{}For every $\mu\in\Mb$, the solution $\um$
in the sense of Stampacchia of the problem
$$
\cases{{\cal A}\um=\mu\ \ & in $\Omega\,$,\crr
	\um=0& on $\partial\Omega\,$,\cr}\eqno(\frm{stampacchia})
$$
is the unique function $\um\in\Luno$ such that
$$
\intl_\Omega \um g\,dx\,=\,\intl_\Omega u^*_g\,d\mu\,,
\quad\hbox{for every }
g\in\Linf\,,
$$
where $u^*_g$ is the solution of
$$
\cases{{\cal A}^*u^*_g=g\quad \hbox{ in }\Hduale\,,\crr
	u^*_g\in\Hunozero\,,\cr}
$$
and ${\cal A}^*$ is the adjoint of ${\cal A}$.
\bigskip

Existence and uniqueness of $u_{\mu}$ are proved in~[\rf{STA}].
Let $T_k(s):=(-k)\vee(s\wedge k)$ be the usual truncation function.
It is easy to prove that
$$
T_k(u_{\mu})\in\Hunozero
\quad\hbox{and}\quad
\intl_\Omega|DT_k(u_{\mu})|^2dx\leq k\,|\mu|(\Omega)\,,
\eqno(\frm{tkho})
$$
for any $k>0$.
These facts imply that $u_{\mu}$ has a quasi continuous representative which is
finite q.e.\ in $\Omega$. If $\mu\in\Inters$, then the solution in 
the sense of Stampacchia coincides with is a 
the usual variational solution in $\Hunozero$.

In the rest of the paper, for every $\mu\in\Mb$ we shall use the 
notation $u_{\mu}$ to indicate the quasi continuous representative of 
the solution of (\rf{stampacchia}), which is uniquely defined quasi 
everywhere in~$\Omega$.

The Green's function $\Gao(x,y)$ relative to the operator ${\cal A}$ in
$\Omega$ is defined as the solution in the sense of Stampacchia of the
equation 
$$
\cases{{\cal A}\Gao(\cdot,y)=\delta_y\ \ & in $\Omega\,$,\crr
	\Gao(\cdot,y)=0& on $\partial\Omega\,$,\cr}
	\eqno(\frm{gaom})
$$
where $\delta_{y}$ is the unit mass concentrated at $y\in\Omega$.
In [\rf{STA}] it is proved that 
$\Gao\colon \Omega\times\Omega\to [0,+\infty]$
is continuous and satisfies the following estimates: for every compact set
$K\subseteq\Omega$ there exist four constants $c_1>0$, $c_2>0$, 
$d_1\geq 0$, and $d_2\geq 0$ ($d_1=d_2=0$ if $N\geq 3$), such that
$$
c_1 G(|x-y|)-d_1\leq \Gao(x,y)\leq c_2 G(|x-y|)+d_2\,,\eqno(\frm{confronto})
$$
for every $x,y\in K$, where $G(|x|)$ is the fundamental solution of
$-\Delta$ in $\Rn$, i.e.,
$$
G(|x|)=\cases{\displaystyle{1\over
		(N-2)\sigma_{N-1}}\,{1\over |x|^{N-2}}\,,
		\quad&if $N>2\,$,\crr
	\displaystyle{1\over 2\pi}\log\big({1\over|x|}\big)\,,
	            &if $N=2\,$,\cr}
\eqno(\frm{fsol})
$$
with $\sigma_{N-1}$ equal to the $(N-1)$-dimensional 
measure of the boundary of
the unit ball in $\Rn$. 
As proved in [\rf{STA}], the
solution of (\rf{stampacchia}) satisfies
$$
\um(x)=\intl_\Omega\Gao(x,y)\,d\mu(y)\,,
\quad \hbox{ for a.e.\ }\,x\in\Omega\,.
\eqno(\frm{verde})
$$

\bigskip

The following notion of solution for obstacle problems with measure
data has been introduced in [\rf{DAL-LEO}].
  
\defin{\thm{disvar}}{}Let $\mu\in\Mb$. We say that a function $u$
is a solution of the
obstacle problem with datum $\mu$ and obstacle $\psi$ (shortly
$OP(\mu,\psi)$) if the following conditions are satisfied:

\item{(a)}$u\in \K$ and there exists $\lambda\in\Mbp$ such
that $u=\um+\ul$ q.e.\ in~$\Omega$;
   
\item{(b)}$u\leq v$ q.e.\ in $\Omega$ for every $v\in\K$
such that  $v=\um+\un$ q.e.\ in $\Omega$, with $\nu\in\Mbp$.
\bigskip

Existence and uniqueness of the solution of $OP(\mu,\psi)$ are proved 
in [\rf{DAL-LEO}], assuming that the obstacle $\psi$ satisfies the 
following natural hypothesis, which replaces (\rf{zeta}):
$$
\hbox{there exists } \rho\in\Mb \hbox{ such that }
u_\rho\geq\psi\qe\,.
$$
In this case we shall say that $\psi$ is {\it $OP$-admissible\/}.
 
The non-negative measure $\lambda$ which appears in condition (a)
of Definition~\rf{disvar} is uniquely determined by the solution $u$ 
and is called the {\it obstacle reaction\/} associated with~$u$.
It is possible to prove that $\lambda$ belongs to $\Mbo$
if the datum
$\mu$ belongs to $\Mbo$ and
$$
\hbox{there exists }  \sigma\in\Mbo \hbox{ such that } 
u_\sigma\geq\psi\qe\,.
$$
When the last condition is satisfied, we shall say that 
$\psi$ is {\it $OP^o$-admissible\/}.
Notice that, if the datum $\mu$ is in $\Mbo$, but the obstacle is
only $OP$-admissible, then the reaction $\lambda$ may
not belong to $\Mbo$. 
For instance, if $\mu=0$ and
$\psi=u_{\delta_{y}}$ for some $y\in\Omega$, then
the solution of $OP(0,\psi)$ is $u_{\delta_{y}}$, and hence
$\lambda=\delta_{y}\not\in\Mbo$.
 
\bigskip

If the obstacle $\psi$ is continuous, or, more in general, quasi 
upper semicontinuous, then the solution of the variational 
inequality (\rf{vi}) must touch the obstacle at all points where it 
is not solution of the equation ${\cal A}u=\mu$.
Indeed, under these assumptions on $\psi$,
the obstacle reaction $\lambda$
of the solution of (\rf{vi}) with  $\mu\in \Hduale$ 
is concentrated on the coincidence set
$\{{x\in\Omega}:{u(x)=\psi(x)}\}$; in other words, $u=\psi$ 
$\lambda$-a.e.\ in $\Omega$. When $\psi$ is continuous, this result 
is well known and can be found in the books mentioned above; the
quasi upper semicontinuous case is discussed, e.g.,
in Section 3 of [\rf{ATT-PIC}].

The same properties are true for the solutions of $OP(\mu,\psi)$ when 
$\mu\in\Mbo$ and $\psi$ is $OP^o$-admissible and 
quasi upper semicontinuous (see [\rf{LEO}]), but they do not hold for 
an arbitrary $\mu\in\Mb$, as shown by the following example, which is 
a variant of an example studied by L.~Orsina and A.~Prignet.

\ex{\thm{delta}}{}Let $\mu\in\Mbp$ be a non-negative measure 
concentrated on a set of capacity zero.
Suppose that there exists a constant $k>0$ such that
$-k\le \psi\le 0$ q.e.\ in $\Omega$.
Let $u=u_{-\mu}+\ul$ be the solution of $OP(-\mu,\psi)$.
We want to show that $u=0$ q.e.\ in $\Omega$ and $\lambda=\mu$
in $\Omega$. 

Taking $\nu=\mu$ in condition (b) of 
Definition~\rf{disvar}, we obtain $u\leq0$ q.e.\ in $\Omega$.
As $u\geq -k$ q.e.\ in $\Omega$, we have $u=T_{k}(u)$ 
q.e.\ in $\Omega$, and hence $u\in\Hunozero$ by (\rf{tkho}). 
This implies that the
measure $-\mu+\lambda$ belongs to $\Inters$, 
which is contained in $\Mbo$. In
other words $\lambda=\mu+\lambda_{0}$, with $\lambda_{0}\in\Mbo$.
Since $\lambda$ is non-negative and $\mu\perp\lambda_{0}$
(recall that $\mu$ is concentrated on a set of capacity zero), the 
measure $\lambda_{0}$ is non-negative. As $u=u_{\lambda_{0}}$, 
by the maximum principle
we have $u\geq 0$ q.e.\ in $\Omega$. Therefore
$u=0$ q.e.\ in $\Omega$ and, consequently, $\lambda=\mu$
in $\Omega$. 

In particular, if $\mu=\delta_{y}$ for some $y\in\Omega$, and $\psi=-k$, 
we have an example of a continuous obstacle for which the solution
$u$ of $OP(-\mu,\psi)$ does not touch $\psi$, although $u$ is 
not the solution of the equation ${\cal A}u=-\mu$, since the obstacle 
reaction is not zero.

\bigskip

In Section~\rf{potenziale} we will show that, when the obstacle is 
controlled from above and from below in an appropriate way 
(see Theorem \rf{assorbimento}), it is possible to
``isolate'' the effect of the singular negative part of the data. 
Namely, the reaction $\lambda$ will be written
as $\lambda=\lambda_0+\mu^-_s$,
where $\lambda_0$ belongs to $\Mbop$. 
Moreover the ``regular part'' $\lambda_0$ is concentrated 
on the coincidence set $\{{x\in\Omega}:{u(x)=\psi(x)}\}$ whenever
$\psi$ is quasi upper semicontinuous, and a complementarity
condition holds (Theorem~\rf{complementarieta}).

The proof of these facts will be based on some new results in 
Potential Theory, which are obtained in the next section.

\goodbreak
\parag{\chp{potenziale}}{Some results in Potential Theory}

We will prove some results concerning the potential of a measure.
The first two lemmas characterize the measures of $\Mbop$ in
terms of the sets where their potentials are infinite. 
The main result of this section is
Lemma~\rf{rapporto} on the behaviour
of the potentials of two mutually singular measures near the points where
both potentials tend to infinity.
It allows us to study the solutions of two equations of the form
(\rf{stampacchia}) corresponding to mutually singular data. In particular
we will compare these solutions near their singular
points (Lemma \rf{rapportoA}).

For every $\mu\in\Mbp$ we consider the potentials $G\mu$ and $\Gao\mu$
defined by
$$
\displaylines{
G\mu(x)=\intl_\Omega G(|x-y|)\,d\mu(y)\,,
\quad\hbox{for } x\in\Rn\,, 
\cr
\Gao\mu(x)=\intl_\Omega \Gao(x,y)\,d\mu(y)\,,
\quad\hbox{for } x\in\Omega\,,
\cr}
$$
where $G$ and $\Gao$ are defined in (\rf{fsol}) and (\rf{gaom}).
Note that $-\Delta G\mu=\mu$ in the sense of distributions in $\Omega$. By
(\rf{verde}) $\Gao\mu$ coincides almost everywhere with the solution 
$\um$ of (\rf{stampacchia}).

\lemma{\thm{cnes}}{Let $\mu\in\Mbp$. Then
$$
\mu\in\Mbop\iff
G\mu<+\infty\quad\mu\hbox{-a.e.\ in }\Omega\,.
$$}
{\vskip-.5cm}
\proof
One implication is easy: by a classical result (see, e.g., Theorem 7.33 in 
[\rf{HEL}]) $G\mu$ is
finite q.e.\ in $\Omega$, and hence $\mu$-a.e.\ in $\Omega$ if
$\mu\in\Mbop$.

Let us prove the converse in the case $N>2$, so that $G\geq 0$.
We start by proving that $\mu^s(\{{x\in\Omega}: {G\mu(x)<+\infty}\})=0$.
For every $t>0$, let $E_t:=\{{x\in\Rn}:{G\mu(x)\leq t}\}$, and let 
$\mu_t$ be the measure defined by
$\mu_t(B):=\mu(B\cap E_t)$ for every Borel set $B\subseteq\Omega$. 
Note that $E_t$  is closed since
$G\mu$ is lower semicontinuous.
As $\mu_t\leq \mu$, we have $G{\mu_t}\leq
G\mu$ (recall that $G\geq 0$). 
In particular
$G{\mu_t}\leq t$ in $E_t$. By the maximum principle (see, e.g.,
Theorem 1.10 in [\rf{LAN}]) we obtain
$G_{\mu_t}\leq t$ in $\Rn$.
Since $G\mu_t$ is superharmonic and bounded, it belongs to ${\rm
H}^1_{\rm loc}(\Rn)$ (see, e.g., Corollary 7.20 in [\rf{HEI-KIL}]). As
$\mu_t=-\Delta G\mu_t$ in
the sense of distributions in $\Omega$, we have $\mu_t\in
\Hduale$, and hence $\mu_t\in\Mbop$.

Let us consider a Borel set $B\subseteq\{{x\in\Omega}:{G\mu(x)<+\infty}\}$ 
with ${\rm cap}(B)=0$.
Then $B$ is the union of the sets $E_t\cap B$, for $t>0$,
and hence 
$$
\mu(B)=\displaystyle\sup_{t\in\R^+}\mu(E_t\cap B)=
\sup_{t\in\R^+}\mu_t(B)=0\,.
$$
Consequently $\mu^s(\{{x\in\Omega}:{G\mu(x)<+\infty}\})=0$.
Therefore, if
$\mu^s$ were not identically zero, it would be 
$\mu^s(\{{x\in\Omega}:{G\mu(x)=+\infty}\})>0$,
and this would contradict the assumption $G\mu<+\infty$ 
$\mu$-a.e.\ in $\Omega$.

The case $N=2$ can be dealt with by adding a suitable constant $c$ to $G$
so that $G+c\geq 0$ in $\Omega$.
The proof is the same with minor modifications, 
among which we point out the use of the maximum
principle for logarithmic potentials 
(see, e.g., Theorem~1.6 in~[\rf{LAN}]).
\endproof

\bigskip

Using (\rf{confronto}) we can now extend Lemma \rf{cnes} to the general
case of the operator~${\cal A}$.

\lemma{\thm{cnesgenerale}}{Let $\mu\in\Mbp$. Then
$$
\mu\in\Mbop\iff
\Gao\mu<+\infty\
\mu\hbox{-a.e.\ in }\Omega\,.
$$
}
{\vskip-.5cm}
\proof
Thanks to (\rf{confronto}) it is easy to prove that for every $x\in\Omega$ 
$$
G\mu(x)<+\infty\iff\Gao\mu(x)<+\infty \,,
\eqno(\frm{sse})
$$
so the thesis follows from Lemma \rf{cnes}.
\endproof

\bigskip

The mean value of an integrable function $f$ on a 
measurable set $B$ with positive measure
is defined by
$$
\ave_{B} f\,dx := {1\over{\rm meas}(B)} \intl_{B} f\,dx\,.
$$
In the next lemma we compare the mean values
of the potentials of two 
mutually singular measures on small balls centered at a point where 
both potentials are infinite.

\lemma{\thm{rapporto}}{Let $\mu,\nu\in\Mbp$,
with $\mu\perp\nu$, and let 
$$
E:=\{x\in\Omega:G\mu(x)=G\nu(x)=+\infty\}\,.
$$
Then
$$
\lim_{r\to0^+}{{\ \ \ \ \displaystyle\ave_{B_r(x)}G\nu\,dy\ \ \ }
\over{\displaystyle\ave_{B_r(x)}G\mu\,dy}}=0\,,
\quad\hbox{ for }\,\mu\hbox{-a.e.\ }x\in E\,.\eqno(\frm{limite1})
$$
}
{\vskip-.5cm}
\proof
Let $R>0$ be such that $\Omega\subseteq B_R(0)$. Observing that
$\Omega\subseteq B_{2R}(x)$ for every $x\in\Omega$, we have
$$
\ave_{B_r(x)}G\nu\,dy=
\intl_{B_{2R}(x)}G_r(|x-z|)\,d\nu(z)\,,
$$
where 
$$
G_r(|x-z|):=\displaystyle\ave_{B_r(x)}G(|y-z|)\,dy\,,
$$
and $\nu$ is
defined for every Borel set $B\subseteq\Rn$ by $\nu(B)=\nu(B\cap\Omega)$. 
As $G(|x|)$ is superharmonic in $\Rn$ and harmonic for $x\neq 0$, we
obtain
$$
G_r(s)\cases{=G(s)\,,& for $s\geq r\,$,\crr
		\leq G(s)\,,& for $s<r\,$,\cr}
$$
and $G_r(s)\nearrow G(s)$ as $r\searrow0$.

It is easy to prove that
$$
\intl_{B_{2R}(x)}G_r(|x-z|)\,d\nu(z) \,=\, G_r(2R)\,\nu(\Omega)
- \intl_0^{2R}G'_r(s)\,\nu(B_s(x))\,ds\,;
\eqno(\frm{spezzo})
$$
the proof can be obtained by using polar coordinates
if $\nu$ is absolutely continuous with respect to the Lebesgue measure, 
and an easy approximation argument
extends the result to the general case.
Note that $\nu(\Omega)<+\infty$ and that $G_r(2R)=G(2R)$ for $r$ small
enough. Since the left hand side of (\rf{spezzo}) tends to $G\nu(x)=+\infty$, 
the last term tends to infinity for every ${x\in E}$.

The same argument can be developed for the denominator, so the limit in
(\rf{limite1}) is equal to
$$
\lim_{r\to0^+}{
{\displaystyle\ \ \ 
\intl_0^{2R}G'_r(s)\,\nu(B_s(x))\,ds\ }
\over{\displaystyle\ \ \ \intl_0^{2R}G'_r(s)\,\mu(B_s(x))\,ds\
}}\,,\eqno(\frm{limite2})
$$
for every ${x\in E}$.
Given $\delta\in (0,2R)$, the integrals between $\delta$ and $2R$
remain
bounded as $r\to0$, so that (\rf{limite2}) is equal to
$$
\lim_{r\to0^+}{
{\displaystyle\ \ \
\intl_0^{\delta}G'_r(s)\,\nu(B_s(x))\,ds\ }
\over{\displaystyle\ \ \ \intl_0^{\delta}G'_r(s)\,\mu(B_s(x))\,ds\ 
}}\,,\eqno(\frm{limite3})
$$
for every ${x\in E}$.
Since $\mu\perp\nu$, by the Besicovitch differentiation theorem (see,
e.g., Chapter 1.6 in
[\rf{EVA-GAR}]), for $\mu$-a.e.\ $x\in\Omega$ we have
$$
\lim_{r\to 0^+}{\nu(B_r(x))\over \mu(B_r(x))}=0\,.
\eqno(\frm{besic})
$$

Let us fix $x\in E$ such that (\rf{besic}) holds. For each $\e>0$ 
there exists $\delta>0$ 
such that 
$$
\nu(B_r(x))<\e\mu(B_r(x)), \quad\hbox{ for all } r\in(0,\delta),
$$
and since $G_r$ is decreasing in $s$, we have
$$
-\intl_0^\delta G'_r(s)\,\nu(B_s(x))\,ds
\leq-\e\intl_0^\delta G'_r(s)\,\mu(B_s(x))\,ds\,.
$$
This shows that the limit in (\rf{limite3}), and hence
in (\rf{limite2}), is
less than or equal to $\e$. Since $\e$ is arbitrary, the limit in
(\rf{limite2}) is zero and the proof is complete.
\endproof

\bigskip
Using (\rf{confronto}) we can extend Lemma~\rf{rapporto}
to the general
case of the operator ${\cal A}$.

\lemma{\thm{rapportoA}}{Let $\mu,\nu\in\Mbp$,
with $\mu\perp\nu$, and let $F$ be the set of all
points $x\in\Omega$ such that
$$
\lim_{r\to 0^+}\ave_{B_r(x)}\um\,dy=
\lim_{r\to 0^+}\ave_{B_r(x)}\un\,dy=
+\infty\,.
\eqno(\frm{limit})
$$
Then  
$$
\lim_{r\to0^+}{{\ \ \ \ \displaystyle\ave_{B_r(x)}\un\,dy\ \ \ }
\over{\displaystyle\ave_{B_r(x)}\um\,dy}}=0\,,
\quad\hbox{ for  }\,\mu\hbox{-a.e.\ }x\in F\,.
\eqno(\frm{B})
$$}
\proof
Let us fix $x\in F$ and $R>0$ such that $B_R(x)\cc\Omega$. 
By (\rf{verde}) we have
$$
\eqalign{
	\ave_{B_r(x)}&\un\,dy=\ave_{B_r(x)}\intl_\Omega
	\Gao(y,z)\,d\nu(z)\,dy\cr
	&=\intl_{\Omega\setminus
	B_R(x)}{\hskip.10cm}\ave_{B_r(x)}\Gao(y,z)\,dy\,d\nu(z)+
	\intl_{B_R(x)}{\hskip.09cm}
	\ave_{B_r(x)}\Gao(y,z)\,dy\,d\nu(z)\,.\cr}
$$
The first term is bounded when $r<{R/ 2}$, so only the second
one is relevant in the limit in (\rf{B}).
The same can be said of the denominator, so that it is enough to study
the quotient
$$
{{\ \ \ \ \displaystyle
\intl_{B_R(x)}{\hskip.15cm}\ave_{B_r(x)}\Gao(y,z)\,dy\,d\nu(z)\ \
\ }
\over{\displaystyle\intl_{B_R(x)}{\hskip.15cm}\ave_{B_r(x)}\Gao(y,z)\,dy\, 
d\mu(z)}}\,.
$$
Thanks to (\rf{confronto}) this is smaller
than or equal to
$$
{{\ \ \ \ \displaystyle
c_2\intl_{B_R(x)}{\hskip.10cm} \ave_{B_r(x)}G(|y-z|)\,dy\,d\nu(z)
+d_1\nu(B_R(x))\
\ \ }
\over{\displaystyle
c_1\intl_{B_R(x)}{\hskip.10cm} \ave_{B_r(x)}G(|y-z|)\,dy\,d\mu(z)
-d_2\mu(B_R(x))}}\,. 
\eqno(\frm{etichetta})
$$
By (\rf{limit}) and (\rf{confronto}), for every $x\in F$ we have
$$
\lim_{r\to0^+}\ave_{B_r(x)}G\mu\,dy=
\lim_{r\to0^+}\ave_{B_r(x)}G\nu\,dy=+\infty\,.
$$
Since $G\mu$ and $G\nu$ are superharmonic, this implies
$G\mu(x)=G\nu(x)=+\infty$ for every $x\in F$. Therefore
Lemma~\rf{rapporto} shows that (\rf{limite1}) holds for
$\mu$-a.e.\ $x\in F$.

Using  once again the fact that 
the integrals over $\Omega\setminus B_R(x)$
remain bounded as $r\to 0^+$, from (\rf{limite1}) 
we obtain that the quotient in (\rf{etichetta})
tends to zero as $r\to 0^+$ for $\mu$-a.e.\ $x\in F$.
\endproof

\lemma{\thm{unico}}{Let $\mu,\,\nu\in\Mb$, let $\lambda\in\Mbo$, and let
$w\in\Huno$. Assume that $\nu\perp\mu^+$ and that
$\um\leq\un+\ul+ w$ a.e.\ in $\Omega$.
Then $\mu^+\in\Mbop$.}
\proof
First of all the measures $\nu$ and $\lambda$ can be assumed to be
non-negative,
replacing them with their positive parts. The function $w$ can be replaced
by $v+h$, where $h$ is the solution of
$$
\cases{
	{\cal A}h=0\quad \hbox{ in }\Hduale\,,\crr
	h-w^+\in\Hunozero\,,\cr}
$$
and  $v=(w-h)^+$. Note that $h$ is a non-negative ${\cal A}$-harmonic
function and $v$ is a non-negative function of $\Hunozero$, and we still
have $\um\leq\un+\ul+v+h$ a.e.\ in $\Omega$.

\medskip
{\it Step 1.} Consider first the case $\um\leq\un$ a.e.\ in $\Omega$. 
Then $u_{\mu^+}\leq\un+u_{\mu^-}$ a.e.\ in $\Omega$,
and $\mu^+\perp(\nu+\mu^-)$.
Let $E$ be the set of all points $x\in\Omega$ such that 
$$
\displaystyle\lim_{r\to 0^+}\ave_{B_r(x)} u_{\mu^+}\,dy=+\infty \,.
$$
Note that $E$ coincides with the set $F$ of Lemma \rf{rapportoA}, 
relative to
the non-negative measures $\mu^+$ and $\nu+\mu^-$. 
Consequently we have
$$
\lim_{r\to0^+}{{\ \ \ \
\displaystyle\ave_{B_r(x)}u_{(\nu+\mu^-)} \,dy\ \ \ }
\over{\displaystyle\ave_{B_r(x)}u_{\mu^+} \,dy}}=0\,,
\quad\hbox{ for }\mu^+\hbox{-a.e.\ }x\in E\,.\eqno(\frm{limiterho})
$$
Since $u_{\mu^+}\leq\un+u_{\mu^-}$ a.e.\ in $\Omega$, the quotient in
(\rf{limiterho}) is greater than or equal to~$1$. Therefore we 
conclude that $\mu^+(E)=0$. As $\Gao\mu^+$ is lower semicontinuous, 
by~(\rf{verde}) we
have $\Gao\mu^+(x)<+\infty$ for $x\in\Omega\setminus E$,  and this
implies $\mu^+\in\Mbop$ by Lemma \rf{cnesgenerale}.
\medskip
{\it Step 2.} Assume that $\um\leq\un+h$ a.e.\ in $\Omega$. Since $h$ is 
${\cal A}$-harmonic, by De Giorgi's theorem it is continuous, hence
$$
\lim_{r\to 0^+}\ave_{B_r(x)}h\,dy=h(x)<+\infty\,,
\quad\hbox{ for every } x\in\Omega\,.
$$
Therefore, if we add this integral  
to the numerator of (\rf{limiterho}),
we can repeat the argument of {\it Step 1\/} and we obtain 
$\mu^+\in\Mbop$ in this
case too.
\medskip
{\it Step 3.} Assume that
$\um\leq\un+h+\ul$ a.e.\ in $\Omega$.
As before we have $u_{\mu^+}\leq u_{(\nu+\mu^-)}+h+\ul$ 
a.e.\ in $\Omega$, with
$\mu^+\perp(\nu+\mu^-)$.
We write now $\mu^+=\mu_1+\mu_2$ in $\Omega$, with $\mu_i\in\Mbp$, 
$\mu_1\ll\lambda$, and $\mu_2\perp\lambda$. Then
$u_{\mu_2}\leq u_{(\lambda+\nu+\mu^-)}+h$  a.e.\ in $\Omega$, 
and $\mu_{2}\perp(\lambda+\nu+\mu^-)$. Therefore we have
$\mu_2=\mu_2^+\in\Mbop$ by {\it Step 2\/}.
As $\mu_1\in\Mbop$, being $\lambda\in\Mbop$,
we conclude that $\mu^+\in\Mbop$.

\medskip
{\it Step 4.} Assume now that  $\um\leq\un+h+\ul+v$ 
a.e.\ in $\Omega$.
Consider the obstacle $\psi_{0}:=\um-\un-h-\ul$, 
which is bounded from above both by $v$
and by $\um$, so that it is both $VI$- and $OP$-admissible.
Then the solution $u_\tau$ of $OP(0,\psi_{0})$ belongs to $\Hunozero$
(see Theorem~5.2  in [\rf{DAL-LEO}]),
hence $\tau\in\Mbp\cap\Hduale\subseteq\Mbop$.
So $\um\leq\un+h+u_{(\lambda+\tau)}$ a.e.\ in $\Omega$, and we 
conclude by means of
{\it Step 3\/}.
\endproof

\cor{\thm{modulo}}{Let $\mu,\,\nu\in\Mb$, let $\lambda\in\Mbo$, and let
$w\in\Huno$. Assume that $\nu\perp\mu$ and that $|\um|\leq\un+\ul+w$
a.e.\ in $\Omega$.
Then $\mu\in\Mbo$.}

\proof
It is enough to apply Lemma~\rf{unico} to $\mu$ and $-\mu$.
\endproof

\parag{\chp{reazione}}{Interaction between obstacles and singular data}

The next theorem is the main result of the paper. We prove that
the component of $\mu^-$ which is singular with respect to the capacity is
completely absorbed by the obstacle reaction $\lambda$, provided 
the obstacle $\psi$ satisfies very weak estimates from above and from 
below.

\th{\thm{assorbimento}}{Let $\mu\in\Mb$ and let $\mu^-_s$ be the part
of $\mu^-$ which is concentrated on a set of capacity zero.
Assume that the obstacle $\psi$ satisfies the estimates
$$
-u_\tau-u_\sigma-w\leq\psi\leq u_\sigma \qe\,,
\eqno(\frm{condizione})
$$
where $w\in\Huno$, $\sigma\in\Mbo$, and $\tau\in\Mb$, with $\tau\perp
\mu^-_s$.
Let $u=\um+\ul$ be the solution of $OP(\mu,\psi)$. Then
$\lambda=\lambda_0+\mu^-_s$ in $\Omega$,
with $\lambda_0\in\Mbop$.}
\proof
It is not restrictive to assume that $\sigma\ge0$ in $\Omega$.
Using the decomposition $\mu^-=\mu^-_a + \mu^-_s$, with 
$\mu^-_a\in\Mbop$,
we can write $u=u_{\mu^+}-u_{\mu^-_a}-u_{\mu^-_s}+\ul$ q.e.\ in
$\Omega$. As 
$\um+u_{(\mu^- +\sigma)}=u_{\mu^+}+u_\sigma\geq\psi$ q.e.\ in 
$\Omega$,
by Definition \rf{disvar} we have $u_{\mu^+}+u_\sigma\geq u$ 
q.e.\ in $\Omega$, hence
$\ul-u_{\mu^-_s}\leq u_\sigma+u_{\mu^-_a}$ q.e.\ in $\Omega$.
By Lemma \rf{unico} this implies $(\lambda-\mu^-_s)^+\in\Mbo$.

On the other hand, $-u_{\mu^-_s}+\ul\geq\psi-u_{\mu^+}+u_{\mu^-_a}$
q.e.\ in $\Omega$, and hence
$u_{(\mu^-_s-\lambda)}\leq u_{\mu^+} +u_\tau +u_\sigma+  w$
q.e.\ in $\Omega$.
Now $(\mu^++\tau)\perp(\mu^-_s-\lambda)^+$, since 
$\mu^+\perp\mu^-$, $\tau\perp\mu^-_s$, and $\lambda\ge0$ in $\Omega$.
So $(\mu^-_s-\lambda)^+\in\Mbop$ by Lemma~\rf{unico}. 

As $(\mu^-_s-\lambda)^-=(\lambda-\mu^-_s)^+\in\Mbop$, we conclude that
$(\mu^-_s-\lambda)\in\Mbo$. 
Therefore $\lambda=\lambda_0+\mu^-_s$, with
$\lambda_0\in\Mbo$. Since
$\lambda\geq 0$ in $\Omega$ and $\lambda_0\perp\mu^-_s$, 
we deduce that $\lambda_0\geq 0$ in~$\Omega$.
\endproof

\rem{\thm{ammiss}}{}Hypothesis (\rf{condizione}) is satisfied,
for instance, when
$\psi$ belongs to $\Huno$ and is $OP$-admissible. Indeed, in this case,
there exists $\rho\in\Mbp$ such that $\psi\leq u_\rho$ q.e.\ in 
$\Omega$.
For any $k\in\R^+$, we have $0\leq\psi^+\wedge k\leq u_\rho\wedge k$. 
Since, by (\rf{tkho}), $u_\rho\wedge k$ 
belongs to $\Hunozero$, so does $\psi^+\wedge k$. As
$$
\intl_\Omega|D(\psi^+\wedge k)|^2dx
\leq\intl_\Omega |D\psi^+|^2dx<+\infty\,,
$$
the function $\psi^+$ is the limit of the increasing sequence
$\psi^+\wedge 
k$, which is bounded in $\Hunozero$. 
This implies that $\psi^+\in\Hunozero$, hence $\psi$ is $VI$-admissible.
Let $u_\sigma$ be the solution of $OP(0,\psi)$. Since $u_\sigma$ is also 
the solution of $VI(0,\psi)$ (see Theorem 5.2 of [\rf{DAL-LEO}]), we have 
$\sigma\in\Inters\subseteq\Mbop$. Then we can take $w=-\psi$ and $\tau=0$
in (\rf{condizione}).

\th{\thm{lostesso}}{Let $\mu\in\Mb$. Assume that the obstacle $\psi$ 
satisfies hypothesis (\rf{condizione}). Let $u$ and $u_{0}$ be the 
solutions of $OP(\mu,\psi)$ and
$OP(\mu^+-\mu^-_a,\psi)$, and let $\lambda$ and $\lambda_{0}$ be the
corresponding obstacle reactions. Then $u=u_{0}$ q.e.\ in $\Omega$ and 
$\lambda=\lambda_{0}+\mu^-_s$ in~$\Omega$. Moreover 
$\lambda_{0}\in\Mbop$.}
\proof
The function $u$ can be written as 
$u_{(\mu^+ - \mu^-_a)} + u_{(-\mu^-_s+\lambda)}$. Since $u\ge \psi$ 
q.e.\ in $\Omega$ and $-\mu^-_s+\lambda\ge 0$ in $\Omega$ by 
Theorem~\rf{assorbimento}, we have $u\ge u_{0}$ q.e.\ in $\Omega$ by 
Definition~\rf{disvar}. Similarly, we have
$u_{0}=u_{(\mu^+ - \mu^-_a)} + u_{\lambda_{0}} = 
u_{\mu} + u_{(\mu^-_s+\lambda_{0})}$
 q.e.\ in $\Omega$. Since $u_{0}\ge \psi$ 
q.e.\ in $\Omega$ and $\mu^-_s+\lambda_{0}\ge 0$ in $\Omega$,
we have $u_{0}\ge u$ q.e.\ in $\Omega$ by 
Definition~\rf{disvar}. Therefore $u= u_{0}$ q.e.\ in $\Omega$ and, 
consequently, $\lambda=\lambda_{0}+\mu^-_s$ in~$\Omega$.
Finally, $\lambda_{0}\in\Mbop$ by Theorem~\rf{assorbimento}.
\endproof

\bigskip

We recall a theorem proved by C. Leone in [\rf{LEO}].

\th{\thm{leone}}{Let $\mu\in\Mbo$ and let $\psi$ be a quasi upper 
semicontinuous $OP^o$-admissible  obstacle. Then the  following facts
are equivalent:

\item{(a)}$u$ is the solution of $OP(\mu,\psi)$ and $\lambda$ is the
corresponding obstacle reaction;

\item{(b)}$\lambda\in\Mbop$, $u=\um+\ul$ q.e.\ in $\Omega$,
$u\geq\psi$ q.e.\ in $\Omega$, and
$u=\psi$ $\lambda$-a.e.\ in $\Omega$.

}

\bigskip

The following theorem extends this result to the case of data in $\Mb$, 
provided the obstacle satisfies~(\rf{condizione}).

\th{\thm{complementarieta}}{Let $\mu\in\Mb$. Assume that
the obstacle $\psi$ is quasi upper semicontinuous and satisfies
hypothesis (\rf{condizione}). 
Then the following facts are equivalent:

\item{(a)}$u$ is the solution of $OP(\mu,\psi)$ and $\lambda$ is the
corresponding obstacle reaction;

\item{(b)}$\lambda=\lambda_0+\mu^-_s$ in $\Omega$, 
with $\lambda_0\in\Mbop$,
$u=\um+\ul$ q.e.\ in $\Omega$,  $u\geq\psi$ q.e.\ in $\Omega$, 
$u=\psi$ $\lambda_0$-a.e.\ in $\Omega$.

}
\proof
{\it Step 1.} First of all we consider the case
$\mu^-\in\Mbop$. 
Observe that $\um+\ul$ is the solution of $OP(\mu,\psi)$ if and only if
$-u_{\mu^-}+\ul$ is the solution of $OP(-\mu^-,\psi-u_{\mu^+})$. By 
Theorem~\rf{leone} this happens if and only if
$\lambda\in\Mbop$, $-u_{\mu^-}+\ul\ge \psi -u_{\mu^+}$ q.e.\ in $\Omega$,
and $-u_{\mu^-}+\ul= \psi -u_{\mu^+}$ $\lambda$-a.e.\ in $\Omega$.
The last two conditions are equivalent to
$u_{\mu}+\ul\ge \psi$ q.e.\ in $\Omega$ and
$u_{\mu}+\ul= \psi$  $\lambda$-a.e.\ in $\Omega$.

\medskip
{\it Step 2.} Let us consider the general case $\mu\in\Mb$.
By Theorem \rf{lostesso} $u$ is the solution of $OP(\mu,\psi)$ and 
$\lambda$ is the corresponding obstacle reaction if and only 
if $u$ is the 
solution of $OP(\mu^+-\mu^-_a,\psi)$ and $\lambda_{0}=\lambda-\mu^-_s$
is the corresponding obstacle reaction.
By {\it Step 1} this happens if and only if $\lambda_{0}\in\Mbop$,
$u=u_{\mu^+} - u_{\mu^-_a} + u_{\lambda_{0}}=\um+\ul$  q.e.\ in $\Omega$,
$u\ge \psi$  q.e.\ in $\Omega$, and 
$u= \psi$ $\lambda_{0}$-a.e.\ in $\Omega$.
\endproof

\intro{References}
\def\interrefspace{\smallskip}  
{\ninepoint\frenchspacing

\item{[\bib{ATT-PIC}]}
Attouch H., Picard C.: Probl\`emes variationnels et
th\'eorie du potentiel non lin\'eaire. {\it Ann. Fac. Sci. Toulouse Math.\/}
{\bf 1} (1979), 89-136.
\interrefspace

\item{[\bib{BAI-CAP}]}
Baiocchi C., Capelo A.: {\it Disequazioni 
variazionali e quasivariazionali. Applicazioni a pro\-ble\-mi di frontiera 
libera.\/} Quaderni dell'Unione Matematica Ita\-lia\-na, Pitagora 
Edi\-tri\-ce, 
Bologna, 1978, translated in {\it Variational and quasivariational 
inequalities.\/} Wiley, New York, 1984.
\interrefspace

\item{[\bib{BOC-CIR1}]}
Boccardo L., Cirmi  G.R.: Nonsmooth unilateral problems.
{\it Nonsmooth optimization: methods and applications
(Erice, 1991), F. Giannessi ed.\/}, 1-10, 
Gordon and Breach, Amsterdam, 1992.
\interrefspace

\item{[\bib{BOC-CIR2}]}
Boccardo L., Cirmi  G.R.:
Existence and uniqueness of solution of unilateral problems with
$L^1$-data. To appear.
\interrefspace

\item{[\bib{BOC-GAL}]}Boccardo  L., Gallou\"et  T.:
Probl\`emes unilat\'eraux avec donn\'ees dans ${\rm L}^1$
{\it C. R. Acad. Sci. Paris S\'er. I Math.\/}  {\bf 311}  (1990), 
617-655.
\interrefspace

\item{[\bib{DAL}]}
Dall'Aglio  P.:
Stability results for solutions of obstacle 
problems with measure data.
In preparation.
\interrefspace

\item{[\bib{DAL-LEO}]}
Dall'Aglio  P., Leone  C.:
Obstacle problems with measure data.
Preprint S.I.S.S.A., Trieste, 1997.
\interrefspace

\item{[\bib{EVA-GAR}]}
Evans  L.C., Gariepy  R.F.:
{\it Measure theory and fine properties of functions.\/}
CRC Press, Boca Raton, 1992.
\interrefspace

\item{[\bib{FUK-SAT-TAN}]}
Fukushima M., Sato K., Taniguchi S.:
On the closable part of pre-Dirichlet forms and the fine supports of
underlying measures.
{\it Osaka J. Math.\/} {\bf 28} (1991), 517-535.
\interrefspace

\item{[\bib{HEI-KIL}]}
Heinonen  J., Kilpel\"ainen  T., Martio  O.: {\it  
Nonlinear potential theory of degenerate elliptic equations.\/} 
Clarendon
Press, Oxford, 1993.
\interrefspace

\item{[\bib{HEL}]}
Helms L.L.:
{\it Introduction to potential theory.\/}
Wiley, New York, 1969.
\interrefspace

\item{[\bib{KIN-STA}]}
Kinderlehrer  D., Stampacchia  G.: {\it
An introduction to variational inequalities and their applications.\/}   
Academic Press, New York, 1980.
\interrefspace

\item{[\bib{LAN}]}
Landkof  N.S.:
{\it Foundations of potential theory.\/} 
Springer Verlag, Berlin, 1972.
\interrefspace

\item{[\bib{LEO}]}
Leone  C.:
Existence and uniqueness of solutions for nonlinear 
obstacle problems with
measure data.
Preprint S.I.S.S.A., Trieste, 1998.
\interrefspace

\item{[\bib{OPP-ROS1}]}
Oppezzi P., Rossi A.M.:
Existence of solutions for unilateral problems with multivalued operators.
{\it J. Convex Anal.\/} {\bf 2} (1995), 241-261.
\interrefspace

\item{[\bib{OPP-ROS2}]}
Oppezzi P., Rossi A.M.:
Esistenza di soluzioni per problemi unilaterali con dato misura o ${\rm L}^1$.
{\it Ricerche Mat.\/}, to appear.
\interrefspace

\item{[\bib{SER}]}
Serrin J.: Pathological solutions of elliptic differential 
equations. {\it Ann.
Scuola Norm. Sup. Pisa Cl. Sci.\/} {\bf 18} (1964), 385-387.
\interrefspace

\item{[\bib{STA}]}
Stampacchia  G.:
Le probl\`eme de Dirichlet pour les \'equations elliptiques du second
ordre \`a coefficients discontinus.
{\it Ann. Inst. Fourier Grenoble\/} {\bf 15} (1965), 189-258.
\interrefspace

\item{[\bib{TRO}]}Troianiello  G.M.: {\it
Elliptic differential equations and obstacle problems.\/}
Plenum Press, New York, 1987.
\par
}

\bye